\documentclass[12pt]{article}

\usepackage[indentafter]{titlesec}
\titleformat{name=\section}{}{\thetitle.}{0.8em}{\centering\scshape}
\titleformat{name=\subsection}[runin]{}{\thetitle.}{0.5em}{\bfseries}[.]
\titleformat{name=\subsubsection}[runin]{}{\thetitle.}{0.5em}{\itshape}[.]
\titleformat{name=\paragraph,numberless}[runin]{}{}{0em}{}[.]
\titlespacing{\paragraph}{0em}{0em}{0.5em}
\titleformat{name=\subparagraph,numberless}[runin]{}{}{0em}{}[.]
\titlespacing{\subparagraph}{0em}{0em}{0.5em}

\usepackage{url}
\makeatletter
\newcommand{\address}[1]{\gdef\@address{#1}}
\newcommand{\email}[1]{\gdef\@email{\url{#1}}}
\newcommand{\sites}[1]{\gdef\@sites{\url{#1}}}
\newcommand{\@endstuff}{\par\vspace{\baselineskip}\noindent\small
\begin{tabular}{@{}l}\scshape\@address\\\textit{E-mail address:} \@email \\ \textit{URL:} \@sites \end{tabular}}
\AtEndDocument{\@endstuff}
\makeatother

\usepackage{titling}
\title{\normalsize\textbf{{\large A} {\large V}IRTUALLY {\large N}ILPOTENT {\large G}ROUP WHOSE \\ {\large G}REEN {\large S}ERIES IS NOT {\large D-F}INITE}}
\author{Corentin Bodart}
\date{\today}

\address{Mathematical Institute, University of Oxford, United Kingdom}
\email{cobodart123@gmail.com}
\sites{https://sites.google.com/view/corentin-bodart}

\usepackage{import}

\usepackage[backend=bibtex, style=numeric, sorting=nyt, 
    doi=false,isbn=false,url=false,eprint=false, maxnames=50]{biblatex}
\usepackage[english]{babel}
\usepackage{amsmath}
\usepackage{amsfonts}
\usepackage{amssymb}
\usepackage{mathrsfs} 

\DeclareMathAlphabet{\funcal}{U}{BOONDOX-cal}{m}{n}
\SetMathAlphabet{\funcal}{bold}{U}{BOONDOX-cal}{b}{n}
\DeclareMathAlphabet{\funbcal} {U}{BOONDOX-cal}{b}{n}

\usepackage{dsfont}
\usepackage{mathtools} 

\usepackage[dvipsnames]{xcolor}

\usepackage{hyperref}
\usepackage{enumitem}
\usepackage[chapter]{easy-todo}         

\usepackage[left=3cm,right=3cm,top=2.8cm,bottom=3cm, marginparsep=3mm]{geometry}
\usepackage{changepage}
\usepackage{adjustbox}

\usepackage[retainorgcmds]{IEEEtrantools}   
\usepackage[font=small, hypcap=false]{caption}
\usepackage{parskip}

\usepackage{tikz-cd}
\usepackage{tikz}
\usetikzlibrary{shapes, shapes.geometric, arrows, arrows.meta, decorations.markings, decorations.pathreplacing, automata}
\pgfdeclarelayer{background} 
\pgfdeclarelayer{foreground}
\pgfsetlayers{background,main,foreground}

\newcommand{\vH}{\funcal{v\!H}}

\newcommand{\say}[1]{``#1"}

\newcommand{\Z}{\mathbb Z}
\newcommand{\Q}{\mathbb Q}

\newcommand{\C}{\mathbb C}
\newcommand{\F}{\mathbb F}

\newcommand{\PP}{\mathbb P} 

\newcommand{\Cay}{\mathcal C{ay}}

\newcommand{\Ac}{\mathcal A}    

\newcommand{\Cc}{\mathcal C}

\newcommand{\Rc}{\mathcal R}

\newcommand{\acts}{\curvearrowright}
\newcommand{\longto}{\longrightarrow}

\newcommand{\la}{\left\langle}
\newcommand{\ra}{\right\rangle}

\newcommand{\abs}[1]{\left| #1 \right|}

\renewcommand{\ge}{\geqslant}
\renewcommand{\le}{\leqslant}

\usepackage{xcolor}

\usepackage{amsthm}
\usepackage{chngcntr}

\theoremstyle{plain}
\newtheorem{thm}{Theorem}
\newtheorem*{thm*}{Theorem}

\newtheorem*{lemma*}{Lemma}
\newtheorem{cor}[thm]{Corollary}
\newtheorem*{cor*}{Corollary}
\newtheorem{prop}[thm]{Proposition}
\newtheorem*{prop*}{Proposition}

\theoremstyle{definition}
\newtheorem{defi}[thm]{Definition}
\newtheorem*{defi*}{Definition}

\newtheorem*{rem*}{Remark}

\begin{document}

\maketitle

\begin{abstract}
	We provide the first example of virtually nilpotent group, with a specific generating set, for which the Green series (sometimes called cogrowth series) is not $D$-finite. The proof relies on an arithmetical miracle, and the study of the subword complexity of a multiplicative sequence coming out of it.
	
	
\end{abstract}
We consider the (re-scaled\footnote{The Green series is usually defined as $G(z)=\sum_{n\ge 0}\mathbb P[X_n=e]\cdot z^n$, where $(X_n)_{n\ge 0}$ is the simple random walk on $\Cay(G,S)$, starting at $X_0=e$. Therefore $\Gamma_{G,S}(z)=G(\abs S \cdot z)$.}) Green series of groups. Given a group $G$ and a finite generating (multi)set $S$, we associate a series $\Gamma_{(G,S)}(z)$ defined as
\[ \Gamma_{(G,S)}(z) = \sum_{\ell\ge 0} c_\ell \cdot z^\ell \in \Z[[z]] \vspace*{-2mm} \]
where $c_\ell=\#\{(s_1,\ldots,s_\ell)\in S^\ell \mid s_1\ldots s_\ell=e_G\}$ is the number of closed paths $e_G\to e_G$ of length $\ell$ in the Cayley graph. In the spirit of the study of other combinatorial sequences, we would like to pin down these series inside the following algebraic hierarchy:
\[ \text{rational} \subset \text{algebraic} \subset \text{diagonal of rational} \subset D\text{-finite} \subset D\text{-algebraic}. \]
Lots of work has been done in this direction. For instance,
\begin{itemize}[leftmargin=8mm]
	\item Kouskov proved that this series is rational if and only if $G$ is finite \cite{Kouskov}.
	\item $\Gamma_{(G,S)}(z)$ is algebraic as soon as $G$ is virtually free. This follows from the Muller-Schupp theorem \cite{Muller_Schupp}, as the word problem is unambiguously context-free in this case. It is an open problem whether the converse holds.
	\item Bishop proved that $\Gamma_{(G,S)}(z)$ is the diagonal of a rational series as soon as $G$ is virtually $F_m\times\Z^n$ \cite{Bishop_cogrowth}. This generalizes previous results for $BS(m,m)$ \cite{BS_cogrowth}.
\end{itemize}
A few results have also been proven in the negative direction:
\begin{itemize}[leftmargin=8mm]
	\item Most proofs use the restricted asymptotics for coefficients of rational/algebraic/$D$-finite growth series. This is the approach used by Kouskov for the result mentioned earlier. The most general result in this direction is that $\Gamma_{(G,S)}(z)$ is not $D$-finite for amenable groups of super-polynomial growth, due to Bell and Mishna \cite{Bell_Mishna}.
	\item Garrabrant and Pak proved that $F_2\times F_2$ (and $SL_4(\Z)$) had a specific generating multiset for which $\Gamma_{(G,S)}(t)$ is not $D$-finite \cite{Garrabrant_Pak}. For $F_2\times F_2$ and symmetric generating set $S$, the asymptotics of $(c_\ell)$ are quite tame, we can show that
	\[ A\cdot \frac{\rho^{2\ell}}{\ell^6} \le c_{2\ell} \le B\cdot \frac{\rho^{2\ell}}{\ell},  \]
	for constants $A,B>0$ and $\rho>1$. (This follows from \cite[Theorem 1.3]{RD_random_walks} and the fact that $F_2\times F_2$ has the Rapid Decay property with exponent $D=3$.) Instead, Garrabrant and Pak develop a new strategy. They construct a generating set so that the sequence $(c_\ell\bmod 4)$ has large \emph{subword complexity}, and prove that this cannot happen for $D$-finite series.
\end{itemize}
For virtually nilpotent groups, the consensus was unclear. 
The known asymptotics perfectly match what is possible among $D$-finite series. For instance, for $G=H_3(\Z)$ with generating set $S=\{x^\pm,y^\pm\}$, it is known that $c_{2\ell}=4^{2\ell}\left(\frac1{2\ell^2}+O\big(\frac1{\ell^3}\big)\right)$ \cite{Gretete}. The first conclusive evidence was given by Pak and Soukup. They encode the existence of solutions to Diophantine equations (which is famously undecidable) inside some decision problem related to Green series. As a corollary, they obtain

\begin{thm*}[Pak-Soukup, \cite{Pak_Soukup}]
	There exists a nilpotent group $G$ such that
	\begin{itemize}[leftmargin=8mm]
		\item Either there exists a finite generating multiset $S$ such that the Green series $\Gamma_{(G,S)}(z)$ is not the diagonal of rational series,
		\item Or at least, there is no algorithm which, given a generating set $S$, computes a representation of the Green series as the diagonal of rational series.
	\end{itemize}	
\end{thm*}
Specifically, they consider the group of $m\times m$ unitriangular matrices $G=UT_m(\Z)$, with $m\approx 10^{86}$. We improve on this result, proving that the first conclusion holds, at least if we allow for virtually nilpotent groups. The group in consideration is
\begin{align*}
	\vH = H_3(\Z)\rtimes C_2
	& = \la x,y,t \;\big|\; [x,[x,y]]=[y,[x,y]]=t^2=e,\, txt = y \ra \\
	& = \la x,t \;\big|\; [x,[x,x^t]]=t^2=e \ra
\end{align*}
introduced in \cite{bishop2020virtHeis} as the first example of group with polynomial geodesic growth which is not virtually abelian. The group $\vH$ is virtually $2$-step nilpotent.
\begin{thm} \label{thm:main}
	The Green series of the virtually nilpotent group $\vH$ with respect to the generating multiset $S=\{x,x^{-1}, t,t,t,t,t,t,t,t\}$ is not D-finite.
\end{thm}
Our proof is most similar to the Garrabrant--Pak argument. We consider the values of some derived sequence modulo a large power of $2$, and compute its subword complexity.
\begin{rem*}
We can promote this result to generating \emph{sets} using $G=\vH\times D_8$ and $T=\{(x^\pm ,0)\}\cup\{(t,g)\mid g\in D_8\}$. This follows from the previous result as
\[ \Gamma_{(\vH,S)}(z)=8\cdot \Gamma_{(G,T)}(z) - 7\cdot \Gamma_{(\la x\ra,\{x^\pm\})}(z), \]
and $\Gamma_{(\la x\ra,\{x^\pm\})}(z)$ is an algebraic series.
\end{rem*}


Observe that we can embed $\vH$ inside $SL_3(\Z)$, sending
\[ x \mapsto \begin{pmatrix} 1 & 1 & 0 \\ & 1 & 1 \\ & & 1 \end{pmatrix},\quad y\mapsto \begin{pmatrix} 1 & 1 & 0 \\ & 1 & \!\!-1 \\ & & 1 \end{pmatrix} \quad\text{and}\quad t\mapsto \begin{pmatrix} -1 & & \\ & \!\!-1 & \\ & & 1 \end{pmatrix}.\]
In particular, using the same trick as Garrabrant and Pak (going from $F_3\times F_{11}$ to any finitely generated overgroup, see the proof of Theorem 1.1 in \cite{Garrabrant_Pak}), we get
\begin{cor}
	$\mathrm{SL}_3(\Z)$ admits a generating multiset such that $\Gamma_{\mathrm{SL}_3(\Z)}(z)$ is not $D$-finite.
\end{cor}
The same result was known in $\mathrm{SL}_m(\Z)$ for $m\ge 4$ \cite{Garrabrant_Pak}, while all Green series of $\mathrm{SL}_2(\Z)$ are algebraic since this group is virtually free.

\section{$D$-finite series}

\begin{defi}
	A series $\sum_{n\ge 0}a_n\cdot z^n \in\Z[[z]]$ is \emph{$D$-finite} if its coefficients are $P$-recursive, i.e., if there exist polynomials $p_0,\ldots,p_k\in\Z[X]$ with $p_0\ne 0$ such that
	\[ p_0(n)\cdot a_n + p_1(n) \cdot a_{n-1} + \ldots + p_k(n) \cdot a_{n-k} = 0 \quad\text{for all } n\ge k. \]
\end{defi}

$D$-finite series are closed under many operations, as shown by Stanley:

\begin{prop}[{\cite[Theorems 2.1, 2.3, 2.7]{STANLEY1980175}}] \label{prop:closure}
	Let $A,\Gamma,\tilde\Gamma\in\Z[[z]]$ be series.
	\begin{enumerate}[leftmargin=8mm, label={\normalfont(\alph*)}]
		\item If $A(z)$ is an algebraic series, then $A(z)$ is $D$-finite.
		\item If $\Gamma(z),\tilde\Gamma(z)$ are $D$-finite, then $c\cdot\Gamma(z)+\tilde c\cdot\tilde\Gamma(z)$ and $\Gamma(z) \cdot \tilde\Gamma(z)$ are $D$-finite.
		\item If $\Gamma(z)$ is D-finite and $A(z)$ is algebraic with $A(0)=0$, then $\Gamma(A(z))$ is D-finite.
	\end{enumerate}
\end{prop}

\medbreak

We add another operation to the list.
\begin{prop} \label{prop:extracted}
	If $\sum_{n\ge 0}a_n\cdot z^n$ is $D$-finite, then the extracted series
	\[ \sum_{n\ge 0} a_{2n}\cdot z^n \quad\text{and}\quad \sum_{n\ge 0}a_{2n+1}\cdot z^n \]
	are $D$-finite.
\end{prop}
\begin{proof} Let $\Gamma(z)=\sum_{n\ge 0}a_n\cdot z^n$. Using the previous Proposition, the series
	\[ \sum_{n\ge 0}a_{2n}\cdot z^{2n} = \frac12\big(\Gamma(z)+\Gamma(-z)\big) \]
	is $D$-finite. Therefore, there exist polynomials $p_0,\ldots,p_{2k}\in\Z[X]$ such that
	\[ p_0(2n)\cdot a_{2n} + p_2(2n)\cdot a_{2n-2} + \ldots + p_{2k}(2n)\cdot a_{2n-2k} = 0 \quad\text{for all }2n\ge 2k.\]
	Taking $q_i(n)=p_{2i}(2n)$ concludes. The proof for the other series is analogous.
\end{proof}

\bigskip
\section{Subword complexity and Multiplicative sequences} \label{sec:complexity}

\begin{defi}
	Given a sequence $(a_n)_{n\ge 0}\in\Ac^{\Z_{\ge 0}}$, its \emph{subword complexity} (or \emph{block complexity}) is the function $p_a\colon\Z_{>0}\to\Z_{> 0}$ defined as
	\[ p_a(n) = \#\big\{(u_1,\ldots,u_n)\in\Ac^n \mid \exists x\ge 0,\; a_{x+i}=u_i \text{ for }i=1,\ldots,n\big\}.\]
\end{defi}

\bigskip

Many \say{algebraically nice} sequences have low complexity. For instance,
\begin{itemize}[leftmargin=8mm]
	\item Eventually periodic sequences (eg.\ coefficients of rational series in $\F_q[[X]]$) are characterized by the property $p_a(n)=O(1)$. Otherwise $p_a(n)\ge n+1$.
	\item Automatic sequences (eg.\ coefficients of algebraic series in $\F_q[[X]]$, or equivalently diagonal of rational series in $\F_q[[X]]$) satisfy $p_a(n)=O(n)$ \cite{frequencies_automatic}.
\end{itemize}
We recall another result in that direction, which will be key in our argument:
\begin{thm}[{\cite[Lemma 4]{Garrabrant_Pak}}] \label{thm:Pak_Garrabrant}
	Let $\sum_{n \ge 0} a_n\cdot X^n\in\Z[[X]]$ be a $D$-finite series. Then the sequence $(a_n\bmod 2)_{n\ge 0}$ has subword complexity $p_a(n)=o(2^n)$.
\end{thm}

\bigskip


In contrast with coefficients of $D$-finite series, we will prove that many \emph{multiplicative} functions have maximal subword complexity. Some of the proof ideas appear \cite[\S 4]{completely_multiplicative_sparse_P}.
\begin{defi}
	A function $f\colon \Z_{>0}\to\C$ is \emph{multiplicative} if
	\[ \forall m,n\in\Z_{>0} \text{ such that } \gcd(m,n)=1,\qquad  f(mn)=f(m)f(n). \]
\end{defi} 
\begin{thm} \label{thm:multiplicative_complexity}
	Let $f\colon \Z_{>0}\to\{\pm1\}$ be a multiplicative function. Suppose that
	\begin{itemize}[leftmargin=8mm]
		\item the set $P_f = \{ p\text{ prime} \mid \exists q= p^m, f(q)=-1\}$ is infinite, and
		\item the set $Q_f = \{ q\text{ prime power} \mid f(q)=-1 \}$ is sparse in the sense $\sum_{q\in Q}\frac 1q<\infty$.
	\end{itemize}
	Then the subword complexity of $\big(f(n)\big)_{n>0}$ is maximal, that is, $p_f(n)=2^n$.
\end{thm}
\begin{rem*}
	Neither assumption can be fully dropped. For instance, if $Q_f$ (hence $P_f$) is finite, then the function is periodic and $p_f(n)=O(1)$.  If we drop the \say{sparseness} condition, some automatic sequences enter the picture, such as
	\[f(n)= \frac{n}{2^{\nu_2(n)}}\pmod 4 \]
	The Liouville function $\lambda(n)$ satisfies the first hypothesis, and $p_\lambda(n)\ge (1+\varepsilon)^n$ is a long-standing open problem, related to Sarnak conjecture on Möbius disjointness.
\end{rem*}

\medskip

\begin{proof}
	Fix $(u_1,\ldots,u_n)\in\{\pm 1\}^n$, we find $x$ such that $f(x+i)=u_i$ for all $1\le i\le n$.
	
	\medbreak
	
	Let $I = \{i\in [\![1,n]\!] \mid f(i)\ne u_i \}$ (the \say{failure set} for $x=0$). For each $i\in I$, we pick a prime power $p_i^{m_i}$ such that $f(p_i^{m_i})=-1$. The first hypothesis ensures we can take all primes $p_i> n$ and distincts. If $i\notin I$, we take as convention $p_i=1$. The Chinese remainder theorem gives infinitely many $x$ satisfying the conditions
	\begin{itemize}[leftmargin=8mm]
		\item For each prime $p\le n$, we take $x \equiv 0 \pmod{p^{m+1}}$, where $m=\lfloor\log_p(n)\rfloor$. 
		\item For each index $i\in I$, we take $x \equiv p_i^{m_i}-i \pmod{p_i^{m_i+1}}$.
	\end{itemize}
	They are all of the form $x=kM+R$ where \vspace*{1mm}
	\[ M=\prod_{p\le n}p^{m+1}\cdot \prod_{i\in I} p_i^{m_i+1}\vspace*{-1mm} \]
	and $0\le R<M$. By construction, we have $x+i=i\cdot p_i^{m_i}\cdot (kM_i+R_i)$ where $kM_i+R_i$ doesn't contain any extra factor $p\le n$ or $p_j$ (with $1\le j\le n$). In particular,
	\[ f(x+i) = f(i)\cdot f(p_i^{m_i})\cdot f(kM_i+R_i) = u_i\cdot f(kM_i+R_i). \]
	We prove that a positive proportion of all $k$ satisfy $q\nmid kM_i+R_i$ for all $q\in Q_f$ and $1\le i\le n$. The only prime factors that still matter come from
	\[ \tilde P_f=\big\{p\in P_f \;\big|\; p> n \text{ and } p\ne p_i \big\}. \]
	For each rank $N$, we partition $\tilde P_f=\tilde P_{f,\le N}\sqcup \tilde P_{f,>N}$. For each $p\in P_{f}$, let $m$ be the smallest integer such that $f(p^m)=-1$
	\[ \frac 1X \#\!\left\{k\le X \;\middle|\; \begin{array}{l}\forall p\in \tilde P_{f,\le N}, \\ \forall i\in[\![1,n]\!], \end{array}\quad p^m\nmid kM_i+R_i\right\} = \prod_{p\in \tilde P_{f,\le N}} \left(1-\frac{n}{p^m}\right) + O_N\left(\frac 1X\right)\]
	(Indeed, the count is exact each time $X$ is a common multiple of the $p^m$ for $p\in \tilde P_{f,\le N}$. The variation in between is accounted by $O_N(\frac1X)$.)
	\begin{align*}
	& \frac1X \cdot \#\! \left\{ k\le X \; \middle|\; \exists p\in \tilde P_{f,>N}, \exists i\in[\![1,n]\!],\quad p^m\mid kM_i+R_i \right\}  \\
	& \hspace*{20mm} \le \frac nX \cdot \#\! \left\{ k\le (M+1)X \;\middle|\; \exists p\in \tilde P_{f,>N},\quad p^m\mid k \right\} \\
	& \hspace*{20mm} \le (M+1)\,n \cdot \sum_{p\in \tilde P_{f,>N}} \frac{1}{p^m}
	\end{align*}
	Using the hypothesis $\sum_{p\in P_f}\frac1{p^m}<\infty$, we have
	\[ \PP_N \coloneqq \prod_{p\in \tilde P_{f,\le N}} \left(1-\frac{n}{p^m}\right) - (M+1)\,n \cdot \sum_{p\in \tilde P_{f,>N}} \frac{1}{p^m} \longto \prod_{p\in\tilde P_f} \left(1-\frac n{p^m}\right) > 0. \]
	Fixing $N$ large enough and letting $X\to\infty$, we get
	\[ \liminf_{X\to\infty} \frac1X \left\{ k\le X\;\middle|\; \forall p\in\tilde P_f, \forall i\in[\![1,n]\!],\quad  p^m\nmid kM_i+R_i \right\} \ge \PP_N > 0 \]
	which concludes.
\end{proof}

\section{Proof of Theorem \ref{thm:main}} \label{sec:main_proof}

\subsection{Model for $H_3(\Z)$}
We recall a path model for the Heisenberg group
\[
H_3(\Z) = \la x,y \mid [x,[x,y]]=[y,[x,y]]=e\ra \simeq
\left\{ \begin{pmatrix}	1 & a & c \\ & 1 & b \\ & & 1 \end{pmatrix} \;\middle|\; a,b,c\in\Z \right\}.
\]
We associate to each word $w$ over $\{x^\pm,y^\pm\}$ a path in $\Z^2$, starting at $(0,0)$. Each letter $x$ (resp.\ $x^{-1}$, $y$, $y^{-1}$) corresponds to a step right (resp.\ left, up, down). We denote the endpoint by $(a,b)$. We concatenate this curve with the path $x^{-a}y^{-b}$, and denote by $c$ the \emph{algebraic area} of the closed curve obtained. Equivalently, $c$ is the sum of the winding numbers of the closed curve around each square of the grid (see Figure \ref{fig:element}). Then the map
\[ w\longmapsto \begin{pmatrix}	1 & a & c \\ & 1 & b \\ & & 1 \end{pmatrix} \]
is an isomorphism. In particular, a word represents the trivial element $e$ if and only if the associated path ends at origin, and has zero algebraic area.
\begin{center}
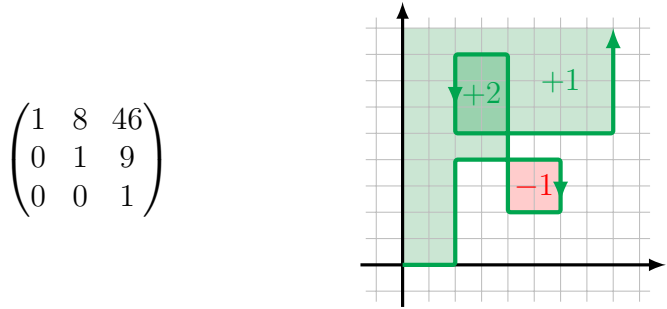

	\begin{tikzpicture}[scale=.7]
		\node at (-6,2) {$\begin{pmatrix} 1 & 8 & \hspace*{-1pt}46 \\ 0 & 1 & \hspace*{-1pt}9 \\ 0 & 0 & \hspace*{-1pt}1 \end{pmatrix}$};
		
		\fill[Green!20] (0,0) -- (1,0) -- (1,2) -- (2,2) -- (2,2.5) -- (4,2.5) -- (4,4.5) -- (0,4.5) -- (0,0);
		\fill[Green!40] (1,2.5) rectangle (2,4);
		\fill[red!20] (2,1) rectangle (3,2);
		\draw[thin, lightgray, step=.5] (-.7,-.7) grid (4.7,4.7);
		\draw[very thick, -latex] (-.8,0) -- (5,0);
		\draw[very thick, -latex] (0,-.8) -- (0,5);
		\draw[Green, ultra thick, -latex, rounded corners=1pt] (0,0) -- (1,0) -- (1,2) -- (3,2) -- (3,1) -- (2,1) -- (2,4) -- (1,4) -- (1,2.5) -- (4,2.5) -- (4,4.5);
		\draw[Green, ultra thick, -latex] (1,3.25) -- (1,3);
		\draw[Green, ultra thick, -latex] (3,1.5) -- (3,1.2);
		\node[Green] at (3,3.5) {$+1$};
		\node[Green] at (1.5,3.25) {$+2$};
		\node[red] at (2.5,1.5) {$-1$};
	\end{tikzpicture}
	\captionsetup{font=footnotesize, margin=30mm}
	\captionof{figure}{The lattice path $x^2y^4x^4y^{-2}x^{-2}y^6x^{-2}y^{-3}x^6y^4$, and the corresponding winding numbers and matrix in $H_3(\Z)$.} \label{fig:element}
\end{center}

\subsection{Reduction to paths without backtracking}

Let us consider the following language (the \say{reduced Word Problem}).
\[  \Rc = \{w\in S^* \mid  \bar w=e,\, \text{no subword }xx^{-1} \text{ or }x^{-1}x \}, \]
and $R(z)=\sum_{\ell\ge 0}r(\ell) \cdot z^\ell$ the associated growth series. Adapting the proof of the Bartholdi--Grigorchuk cogrowth formula \cite[Corollary 2.6]{Bartholdi}, we get
\[ \frac{R(z)}{1-z^2} = \frac{\Gamma\left( \frac z{1+z^2}\right)}{1+z^2}. \]
(As we only remove \say{bumps} $xx^{-1}$ and $x^{-1}x$, we should take $d=2$ in the formula.) It follows from Proposition \ref{prop:closure}(c) that $R(z)$ is $D$-finite if and only if $\Gamma(z)$ is $D$-finite.


\subsection{Counting paths with few $t$'s} 

We decompose $\Rc$ into three disjoint sets:
\begin{align*}
	\Rc_1 & = \left\{ w\in \Rc \;\big|\; \text{at most four }t\text{, or six }t\text{ including two consecutive} \right\}, \\
	\Rc_2 & = \left\{ w\in \Rc \;\big|\; \text{exactly six }t \text{ and no subword } tt \right\}, \\
	\Rc_3 & = \left\{ w\in \Rc \;\big|\; \text{at least eight }t\right\}.
\end{align*}
(1) Observe that $\Rc_1 = WP\big(F_2\rtimes C_2,\{x,x^{-1},8\cdot t\}\big)\cap \Cc$ where $\Cc$ is a rational language encoding the fact that our words are reduced and the condition on the $t$'s. Indeed, if
\[ x^{n_0}\,t\,x^{n_1}\,t\,x^{n_2}\,t\,x^{n_3}\,t\,x^{n_4} = x^{n_0}y^{n_1}x^{n_2}y^{n_3}x^{n_4} \]
is trivial in $\vH$, then it is also trivial in $F_2\rtimes C_2$. Using the easy direction of Müller--Schupp's theorem, we conclude that $\Rc_1$ is unambiguously context-free.

\bigskip 

(2) Paths of length $2\ell+6$ in $\Rc_2$ come in two types and four orientations (Figure \ref{fig:orientations}).

\begin{center}
	\begin{tikzpicture}[xscale=.7, yscale=-.7]
		\draw[step=.5, dotted] (-2.7,-2.2) grid (1.7,3.7);
		\draw[very thick, rounded corners=1] (.97,0) -- (1,0) -- (1,3) -- (0,3) -- (0,2.97)
		(-1.97,0) -- (-2,0) -- (-2,-1.5) -- (0,-1.5) -- (0,-1.47);
		\draw[line width=1.5pt, red] (.07,0) -- (.97,0);
		\draw[line width=1.5pt, orange] (0,0) -- (0,2.97);
		\draw[line width=1.5pt, blue] (-.07,0) -- (-1.97,0);
		\draw[line width=1.5pt, cyan] (0,0) -- (0,-1.47);
		\draw[ultra thick, red, -latex] (.5,0) -- (.75,0);
		\draw[ultra thick, cyan, -latex] (0,-.75) -- (0,-1);
		
		\node[red] at (.75,-.26) {$a$};
		\node[orange] at (-.25,1.55) {$b$};
		\node[blue] at (-1,.25) {$c$};
		\node[cyan] at (0.25,-.93) {$d$};
		
		\node[circle, draw=black, fill=white, inner sep=1.5pt] at (-.5, 3.6) {\scriptsize$\nearrow$};
		
		\node[circle, draw=black, fill=Green, inner sep=1.4pt] at (-1.5,0) {};
	\end{tikzpicture} \hspace*{4mm}
	\begin{tikzpicture}[scale=.7]
		\draw[step=.5, dotted] (-2.7,-2.2) grid (1.7,3.7);
		\draw[very thick, rounded corners=1] (.97,0) -- (1,0) -- (1,3) -- (0,3) -- (0,2.97)
		(-1.97,0) -- (-2,0) -- (-2,-1.5) -- (0,-1.5) -- (0,-1.47);
		\draw[line width=1.5pt, red] (.07,0) -- (.97,0);
		\draw[line width=1.5pt, orange] (0,0) -- (0,2.97);
		\draw[line width=1.5pt, blue] (-.07,0) -- (-1.97,0);
		\draw[line width=1.5pt, cyan] (0,0) -- (0,-1.47);
		\draw[ultra thick, blue, -latex] (-1.05,0) -- (-1.25,0);
		\draw[ultra thick, orange, -latex] (0,1.5) -- (0,1.7);
		
		\node[red] at (.5,-.26) {$a$};
		\node[orange] at (-.25,1.74) {$b$};
		\node[blue] at (-1.25,.24) {$c$};
		\node[cyan] at (0.25,-.77) {$d$};
		
		\node[circle, draw=black, fill=white, inner sep=1.5pt] at (-.5, -2.1) {\scriptsize$\nwarrow$};
		
		\node[circle, draw=black, fill=Green, inner sep=1.4pt] at (0,3) {};
	\end{tikzpicture} \hspace*{4mm}
	\begin{tikzpicture}[xscale=.7, yscale=-.7]
		\draw[step=.5, dotted] (-3.4,-2.2) grid (1.7,3.7);
		\draw[very thick, rounded corners=1] (.97,0) -- (1,0) -- (1,3) -- (0,3) -- (0,2.97)
		 (-2,-1) -- (-2,-1.5) -- (0,-1.5) -- (0,-1.47);
		\draw[very thick, rounded corners=3] (-3,0) -- (-2,0) -- (-2,-1);
		\draw[line width=1.5pt, red] (.07,0) -- (.97,0);
		\draw[line width=1.5pt, orange] (0,0) -- (0,2.97);
		\draw[line width=1.5pt, blue] (-.07,0) -- (-2.05,0);
		\draw[line width=1.5pt, cyan] (0,0) -- (0,-1.47);
		\draw[ultra thick, blue, -latex] (-1.05,0) -- (-1.25,0);
		\draw[ultra thick, orange, -latex] (0,1.5) -- (0,1.7);
		
		\node[red] at (.5,-.26) {$a$};
		\node[orange] at (-.25,1.55) {$b$};
		\node[blue] at (-1.25,.26) {$c$};
		\node[cyan] at (0.25,-.74) {$d$};
		
		\node[circle, draw=black, fill=white, inner sep=1.5pt] at (-.5, 3.6) {\scriptsize$\swarrow$};
		
		\node[circle, draw=black, fill=Green, inner sep=1.2pt] at (-3,0) {};
	\end{tikzpicture} \hspace*{4mm}
	\begin{tikzpicture}[scale=.7]
		\draw[step=.5, dotted] (-2.7,-2.2) grid (1.7,3.7);
		\draw[very thick, rounded corners=1] (.97,0) -- (1,0) -- (1,3) -- (0,3) -- (0,2.97)
		(-1.97,0) -- (-2,0) -- (-2,-1.5) -- (0,-1.5) -- (1,-1.5);
		\draw[very thick, rounded corners=3] (0,0) -- (0,-1.5) -- (1,-1.5);
		\draw[line width=1.5pt, red] (.07,0) -- (.97,0);
		\draw[line width=1.5pt, orange] (0,0) -- (0,2.97);
		\draw[line width=1.5pt, blue] (-.07,0) -- (-1.97,0);
		\draw[line width=1.5pt, cyan, rounded corners=2] (0,0) -- (0,-1.4);
		\draw[ultra thick, red, -latex] (.5,0) -- (.75,0);
		\draw[ultra thick, cyan, -latex] (0,-.75) -- (0,-1);
		
		\node[red] at (.75,-.26) {$a$};
		\node[orange] at (-.25,1.55) {$b$};
		\node[blue] at (-1,.25) {$c$};
		\node[cyan] at (0.25,-.81) {$d$};
		
		\node[circle, draw=black, fill=white, inner sep=1.5pt] at (-.5, -2.1) {\scriptsize$\searrow$};
		
		\node[circle, draw=black, fill=Green, inner sep=1.2pt] at (1,-1.5) {};
	\end{tikzpicture}
	\captionsetup{margin=4mm, font=footnotesize}
	
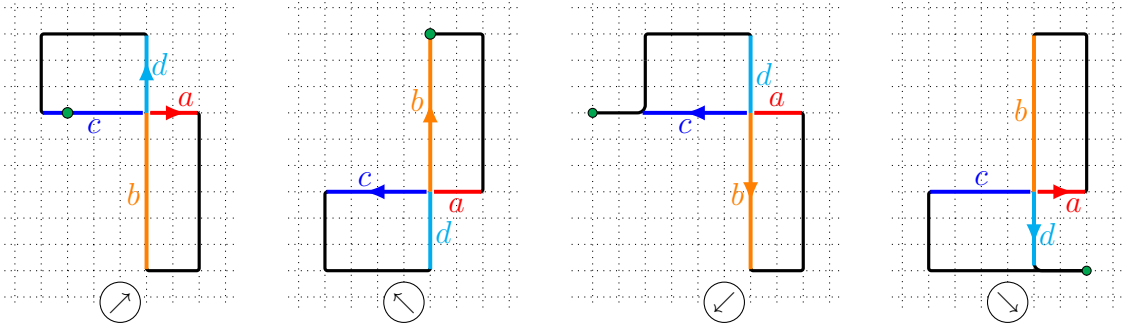
\captionof{figure}{Two paths of each type, and all four orientations. For the first type, the perimeter should be $2\ell$. For the second type, the perimeter of the \say{main shape} (i.e., without the \say{tail}) should be $2k<2\ell$. The first picture corresponds to the word $w=x^5tx^{-6}tx^{-2}tx^9tx^{-4}tx^{-3}tx$.} \label{fig:orientations}
\end{center}
Shapes of perimeter $2n$ are parametrized by solutions $(a,b,c,d)\in\Z_{>0}^4$ to the system
\[ \begin{cases} a+b+c+d=n, \\ ab=cd. \end{cases}\]
We denote the set of solution by $S_n$.

For paths of the first type, we have $2(a+c)+3$ starting points (anywhere along an horizontal segment). For the second type, we have exactly $6$ starting points (the \say{tail} can be attached at any corner, its length $\ell-k$ is fixed). This leads to
\[ r_2(2\ell+6) = 8^6 \cdot \left(\sum_{(a,b,c,d)\in S_\ell} 4 \cdot \big(2(a+c)+3\big) + \sum_{k<\ell} 4 \cdot 6 \cdot \abs{S_k}\right) \]

(3) Finally $8^8\mid r_3(\ell)$ as we can choose any of the $8$ copies for each instance of $t$.

\subsection{Counting solutions to a Diophantine equation} We compute $\abs{S_n}\pmod 8$. A key observation is that, in addition to the four orientations used earlier, there is an extra symmetry. Specifically $D_8\acts S_n$ generated by the involutions
\[ \sigma(a,b,c,d)=(a,b,d,c) \quad\text{and}\quad \tau(a,b,c,d)=(c,d,a,b). \]
Using the Orbit-Stabilizer formula, we get
\[\begin{array}{rclcl}
	\abs{S_n}= & & 1\cdot  \#\{(a,a,a,a)\} & + &  4\cdot \#\{(a,b,a,b):a<b\} \\
	& + & 4\cdot \#\{(a,b,c,c):a<b\} & + & 8\cdot \#\{(a,b,c,d):a<c<d<b\} \\[2mm]
	= & & \mathds 1_{\{4\mid n\}} & + & 4\cdot \mathds 1_{\{2\mid n\}} \cdot \left\lfloor\frac{n-1}4\right\rfloor \\
	& + & 4 \cdot \#\{(a,b,c,c):a<b\} & + & 8 \cdot \text{an integer} 
 \end{array}\]
For each $n\in\mathbb Z_{>0}$, let's compute the number of solutions $(a,b,c)\in\Z_{>0}$ to the system
\[ \begin{cases} a+b+2c = n, \\ \hspace*{7mm} ab \hspace*{7mm} = c^2, \\ \hspace*{4mm} a<b . \end{cases}\]
Let $d=\gcd(a,b)$. We can write $a=dX$ and $b=dY$ with $X<Y$ coprime integers. As $d^2XY=ab=c^2$, we conclude that both $X$ and $Y$ are perfect squares, more precisely $a=dx^2$, $b=dy^2$ and $c=dxy$. Now the first equation becomes
\[ d(x+y)^2 = n \]
Reciprocally, for each integer $z\ge 3$ such that $z^2\mid n$, we have $\frac12\varphi(z)$ choices for $x<y$ such that $x+y=z$ and $\gcd(x,y)=1$, where $\varphi$ is the Euler's totient function. Using Gauss formula $\sum_{d\mid m}\varphi(d)=m$, we conclude that
\[ \#\{(a,b,c,c)\in S_n:a<b\} = \frac12\sum_{\substack{z^2\mid n \\ z\ge 3}}\varphi(z) = \frac12\sum_{\substack{z\mid m \\ z\ge 3}} \varphi(z) = \frac{m(n)-1-\mathds 1_{\{4\mid n\}}}2, \]
where $m(n)=\prod_p p^{\lfloor v_p(n)/2\rfloor}$ is the largest integer such that $m^2\mid n$.

\subsection{Conclusion} Let's put everything together for $\ell=2j+1$ odd:
\begin{align*}
	r(4j+8) & \equiv r_1(4j+8) + 8^6\cdot 4\cdot \left( 4\cdot\frac{m(2j+1)-1}2 + \sum_{k<2j+1} 6 \cdot \mathds 1_{\{4\mid k\}}\right) \pmod{8^6\cdot 4\cdot 8} \\[2mm]
 & \equiv r_1(2\ell+6) +
  2^{21} \cdot \left( m(2j+1)-1 + 3 \cdot \left\lfloor\frac j2\right\rfloor \right)
\end{align*}
hence
\[ \frac1{2^{22}} \left(r(4j+8) - r_1(4j+8) - 2^{21}\cdot 3\cdot \left\lfloor\frac j2\right\rfloor \right) \equiv \frac{m(2j+1)-1}2 \pmod 2. \]
Let
\[ S(z) \coloneqq \sum_{j=0}^\infty \frac1{2^{22}} \left(r(4j+8) - r_1(4j+8) - 2^{21}\cdot 3\cdot \left\lfloor\frac j2\right\rfloor \right) \cdot z^j \in\Z[[z]]. \]
Observe that
\begin{align*}
	& \frac{m(2j+1)-1}2 \equiv 0 \pmod 2 \iff f(2j+1)\equiv +1 \pmod 4 \quad\text{and} \\
	& \frac{m(2j+1)-1}2 \equiv 1 \pmod 2 \iff f(2j+1)\equiv -1 \pmod 4,
\end{align*}
where $f\colon\Z_{>0}\to\{\pm1\}$ is the multiplicative function defined as
\[ f(n)=m\!\left(\frac n{2^{\nu_2(n)}}\right)\pmod 4, \] 
which satisfies the hypothesis of Theorem \ref{thm:multiplicative_complexity}. We conclude that the subword complexity of the coefficients of $S(z)$ modulo $2$ is $p(n)=2^n$, hence $S(z)$ cannot be $D$-finite by Garranbrant--Pak's Theorem \ref{thm:Pak_Garrabrant}. As the generating series of $r_1(4j+8)$ and $\left\lfloor\frac j2\right\rfloor$ are algebraic, it follows that $R(z)$ and $\Gamma(z)$ are not $D$-finite (Proposition \ref{prop:closure} and \ref{prop:extracted}). \qed
\section{Final remarks}

\subsection{Alternate arguments} It is possible to bypass Section \ref{sec:complexity} and get the weaker conclusion that the Green series cannot be written as a diagonal of rational series.

We first repeat the argument of Section \ref{sec:main_proof}: if $\Gamma(z)$ is the diagonal of rational series, then $S(z)$ is too. Using \cite[Theorem 5.2]{diagonal_implies_algebraic} and \cite[Théorème 1]{algebraic_implies_automatic}, we get that the sequence $\big(m(n)\bmod 4\big)_{n\ge 1}$ is $2$-automatic, and therefore its multiplicative cousin 
\[
f(n) = m\!\left(\frac{n}{2 ^{\nu_2(n)}}\right) \pmod 4
\]
(with values in $\{\pm1\}$) is $2$-automatic too. However, sequences that are both automatic and multiplicative are classified \cite{automatic_multiplicative}, and $f$ does not appear on the list.

\bigbreak

Another tempting argument was to use asymptotic frequencies. Unfortunately, this reduces to a well-known open problem. Recall that squarefree numbers have asymptotic density $\frac6{\pi ^2}$. It follows that naturals such that $m(n)=k$ have density $\frac1{k^2}\cdot\frac6{\pi^2}$, hence
\[ \lim_{X\to\infty}\frac{\#\{n\le X \mid m(n)\equiv 1 \bmod 4 \}}X = \frac6{\pi^2}\sum_{k\ge 0}\frac1{(4k+1)^2} = \frac3{\pi ^2}\left(\frac{\pi^2}8+G\right) = \frac38 + \frac{3G}{\pi^2}\]
where $\frac{\pi^2}8=\sum_{\ell\ge 0}\frac1{(2\ell+1)^2}$ and $G=L(2,\chi_{-4})=\sum_{\ell\ge 0}\frac{(-1)^\ell}{(2\ell+1)^2}$ is the Catalan constant. It is widely believed that $\pi^2$ and $G$ are $\Q$-linearly independent (\cite{zeta} for recent progress), hence the sequence $(m(n)\bmod 4)_{n\ge 1}$ cannot be $2$-automatic \cite[Theorem 6]{frequencies_automatic}.

\subsection{The discrete Heisenberg group} The most tempting question is to extend the result to the discrete Heisenberg group
\[ H_3(\Z) = \la x,y\mid [x,[x,y]]=[y,[x,y]] = e \ra \]
itself. Part of the motivation is that any proof looking at the coefficients $\pmod{p^m}$ would then pass to any group containing $H_3(\Z)$, for instance any virtually nilpotent group which is not virtually abelian (for a well-chosen generating multiset).

\medbreak

However, the small trick of adding multiple copies of a generator doesn't seem to work. For instance, if we take $S=\{x,y,y,z\}^{\pm}$ as our generating set and look modulo $2^{2K}$, we have to count closed paths that stay within a $K$-neighborhood of the abelian subgroup $\la x,z\ra$, hence the associated series should be the diagonal of a rational series.

This means that we need to find extra symmetries on the entire set of closed paths (and not just $\Rc_2$ in our argument), for instance find some $2$-group action. Hopefully, once filtering by increasing orbit size, the first terms would be provably good (i.e.\ $D$-finite), and we can find a provably bad term before getting stuck on unprovably bad ones.



\bigskip

\textbf{Acknowledgment.} I would like to thank to Tatiana Nagnibeda for pointing me to the article of Pak and Soukup in the first place, as well as Alex Bishop for useful discussions. The author acknowledges support of the Swiss NSF grant 200020-200400.



\AtNextBibliography{\small}
\printbibliography
\end{document}